\magnification=\magstep1

\documentstyle{amsppt}

\NoRunningHeads \font\tyt=cmbx12

\def\Om{\Omega }

\def\MA{Monge-Amp\`ere }

\def\cl{\centerline }

\def\fii{\varphi }

\def\we{\wedge }

\def\om{\omega }
\def\db{\overline{\partial } }
\def\dd{\partial }
\def\ep{\epsilon }
\def\ga{\gamma }

\def\de{\delta }
\def\ze{\zeta }
\def\te{\theta }

\def\sk{\nopagebreak\vskip 3mm}
\def\skk{\nopagebreak\vskip 5mm\nopagebreak}


\topmatter

\cl{\tyt A uniform $L^\infty$ estimate}

\sk

\cl{\tyt for complex Monge-Amp\`ere equations}

\author S\l awomir Ko\l odziej \ and \  Gang Tian \endauthor


\thanks First author partially supported by EU grant
MTKD-CT-2006-042360 and Polish   grant 189/6. PRUE/2007/7
\endthanks

\endtopmatter

\skk\skk

\bigskip
\bigskip
\bigskip

\document

\proclaim{\bf 0. Introduction} \endproclaim \rm Inspired by works
on the K\"ahler-Ricci flow [TZ], [ST1] and on regularity of
complex Monge-Ampere equations [K1], the second named author made
the following conjecture: Let $\pi : X\to Y$ be a holomorphic
mapping between compact K\"ahler manifolds with $dim X=n\geq m=dim
Y$. Let $\om _X $ be a K\"ahler metric on $X$ and $\om _Y $ be a
K\"ahler metric on $Y$. Suppose that $F\in L^{\infty} (M)$
\footnote{It may be sufficient to assume that $F\in L^p(M)$ for
some $p> 1$.} and $u_t \in PSH(\om _t )$ (which means that the
current $\om _t +dd^c u_t$ is weakly positive) be a solution of
$$\aligned
(\om _t +dd^c u_t )^n= c_t t^{n-m} F\om _X^n , \ \ \ \max _X u_t
=0,  \endaligned\tag{0.1}
$$
where $ \om =\pi ^{\ast} \om _Y , t\in (0,1),$ and $ c_t$ is
defined by

$$\frac{\int _X \om _t ^n }{t^{n-m}\int _X F\om _X^n }.
$$
\rm Then $u_t$ are uniformly bounded. Note that $\om$ is a weakly
positive $(1,1)$ form on $X$ and $ \om _t =\om +t\om _X$ are
K\"ahler forms for any $t>0.$ The main purpose of this note is to
affirm this conjecture under certain technical conditions on the
singular set of $\pi$. We do not know if these conditions are
satisfied for any holomorphic fibration as above. However, we will
check that these conditions are indeed true for many holomorphic
fibarions.

Now let us state our conditions. Assume that for an analytic sets
$S\subset X$ and $A\subset Y$ the form $\om$ is non degenerate
away from $S\subset \pi ^{-1} (A)$. Consider the stratification of
$A$, into $A_1 , A_2 , ... , A_{N_1}$ so that $A_0$ is the set of
regular points of $A$, $A_1$ is the regular part of $A\setminus
A_0$ and so on. Analogously we stratify $S$ into $S_1 , S_2 , ...
, S_{N_2}.$ Let  $g_{q,j,k}$ be a fixed local system of generators
of the ideal sheaf of $\bar{A}_q$ in coordinate patch $V_j$ and
let $\te _j $ be smooth functions on $Y$ such that $supp\,  \te _j
\subset V_j $, $0\leq \te _j \leq 1$ and the interiors $V'_j$ of
$\{ \te _j =1 \}$ cover $Y$. In the same way we define local
generators $h_{q,j,k}$ for ${S}_q$ in $U_j \subset X$ and smooth
functions $\ze _j$ which are equal to $1$ on $U'_j .$

Assume that for some $a\in (0,1) , c_0 >0$ and any fixed $c\in
(0,c_0 )$ we have

$$\aligned & \pi ^{\ast } [ \frac{1}{2} \om _Y + ica^2 \sum _{q,k}
|g_{q,j_0 ,k} |^{2(a-1)}
\partial g_{q,j_0 ,k} \we \bar{\partial }\bar{g}_{q,j_0 ,k} ] ^{m}
 \\ &  \we  [ \frac{1}{2} \om _X + ica^2 \sum _{q,k} |h_{q,j_1 ,k}
|^{2(a-1)}
\partial h_{q,j_1 ,k} \we \bar{\partial }\bar{h}_{q,j_1 ,k}  ]^{n-m}\\
&  \geq H(\ep ) \om _X ^n , \endaligned\tag{0.2} $$
on the set
$U'_{j_1}\cap \pi ^{-1} (V'_{j_0 } )\cap S ^{\ep } ,$ where $S
^{\ep }=\{ \max _{j,k} \ze _j |h_{1,j ,k}| < \ep \}$,  with some
function $H$ satisfying $\lim _{\ep \to 0} H(\ep )=\infty .$

Note that a bit more explicit inequality (on the same sets, with
the same $a,c_0 ,H$ as in \thetag{0.2}):
$$
\aligned & \pi ^{\ast } [ i \sum _{q,k} |g_{q,j_0 ,k} |^{2(a-1)}
\partial g_{q,j_0 ,k} \we \bar{\partial }\bar{g}_{q,j_0 ,k} +\om _Y ]\we \om
^{m-1} \\ &  \we  [ i \sum _{q,k} |h_{q,j_1 ,k} |^{2(a-1)}
\partial h_{q,j_1 ,k} \we \bar{\partial }\bar{h}_{q,j_1 ,k} +\om _X ]\we \om
_X ^{n-m-1}\\
&  \geq H(\ep ) \om _X ^n , \endaligned\tag{0.3} $$ entails
\thetag{0.2}.

The following theorem confirms the above conjecture under the
assumption \thetag{0.2}.

\proclaim{Theorem 1} Suppose that $F\in L^{\infty} (M)$ and the
fibration $\pi:X\mapsto Y$ satisfies \thetag{0.2}, then the
solutions $u_t$ of \thetag{0.1} are uniformly bounded.
\endproclaim

For its application, we consider the expanding K\"ahler-Ricci flow
on $X$
$$ \frac{\partial \omega(t,\cdot)}{\partial t} = -
 Ric(\omega(t,\cdot))- \omega(t,\cdot),
~~~~\omega(0,\cdot)=\omega_0. \tag{0.4}
$$
Let $\pi: X\mapsto Y$ be a holomorphic fibration such that
$c_1(X)=\pi^*\omega_Y $. It follows that smooth fibers of $\pi$
are smooth Calabi-Yau manifolds. Let $Y_0$ be the dense-open
subset of $Y$ over which fibers of $\pi$ are smooth, then there is
an induced holomorphic map $f$ from $Y_0$ into a moduli space of
Calabi-Yau manifolds. This map assigns each $y\in Y_0$ to the
Calabi-Yau manifold $\pi^{-1}(y)$. By [TZ], for any $\omega_0$,
then \thetag{0.4} has a global solution $\omega(t,\cdot)$.
Moreover, there is a smooth family of functions $\varphi(t,\cdot)$
satisfying
$$\omega(t,\cdot)=e^{-t} \omega_0 + (1-e^{-t})\pi^*\omega_Y + i\partial\overline{\partial}
\varphi(t,\cdot).$$

It was expected (cf. [ST2]) that $\omega(t,\cdot)$ converges to a
generalized K\"ahler-Einstein metric on $Y$ in a suitable sense as
$t$ tends to $\infty$. Combining estimates in [ST1] and the above
theorem, we have

\newpage

\proclaim{Theorem 2} Suppose that the fibration $\pi: X\mapsto Y$
has no multiple fibers and satisfies \thetag{0.2}, then there is a
positive current $\omega_\infty$ on $Y$ with properties:

1. $\omega_\infty=\omega_Y + i\partial\overline{\partial}
\varphi_\infty$ for some bounded function $\varphi_\infty$ on $Y$;

2. $\varphi_\infty$ is smooth on $Y_0$;

3. $\varphi(t,\cdot)$ converges to $\pi^*\varphi_\infty$ in
$C^{1,1}$-topology on  any compact subsets contained in
$\pi^{-1}(Y_0)$;

4. On $Y_0$, $\omega_\infty$ satisfies the equation for
generalized K\"ahler-Einstein metrics

$$\text{ Ric }(\omega_{\infty }) = - \omega_{\infty } + f^*\omega_{WP},$$

where $\omega_{WP}$ denotes the Weil-Petersson metric on the
moduli of Calabi-Yau manifolds.
\endproclaim

This theorem can be proved by those estimates developed in [ST1]
without using Theorem 1 if the base is 1-dimensional or $X$ is of
dimension $2$. The rest of this note is organized as follows: In
Section two, we give a proof of Theorem 1. In Section three, we
prove Theorem 2 using Theorem 1 and results in [SZ1] and [SZ2]. In
last section, we verify the assumptions \thetag{0.2} for certain
holomorphic fibrations, including any fibration with 1-dimensional
base and generic 3-folds with 2-dimensional base and tori as
fibers.

\newpage

\proclaim{\bf 1. Proof of Theorem 1} \endproclaim \rm

\demo{Proof} \bf In what follows $C$ denotes different positive
 constants independent of $t$. \rm
 Note that
 $$
 \lim _{t\to 0} c_t = \binom{n}{m}\frac{\int _X \om ^m \we\om _X ^{n-m} }{\int _X F\om _X^n }>0 .
$$

\proclaim{Lemma 1}  For $(Y,\om _Y ), \te _j , V'_j$ as above
there exists $c_1 >0$ such that for all $c\in (0,c_1 )$, $a\in
(0,1)$, and
$$
\psi _a = \sum _{q,j,k} \te _j ^2  |g_{q,j,k}|^{2a} ,
$$
 we have
$$
 dd^c c\psi _a \geq -\frac{1}{2}\om _Y +  ca^2 i\sum _{q,k}   |g_{q,j,k}|^{2(a-1)}
\partial g_{q,j ,k} \we \bar{\partial }\bar{g}_{q,j ,k} \tag{1.1}
$$
on the set $V'_j .$
\endproclaim

\demo{Proof} By computation (comp. [DP, Lemma 2.1])
$$\aligned
i\dd\db \psi _a &= 2i\sum _{q,j,k}
 |g_{q,j,k}|^{2a} ( \te _j \dd\db\te _j  -\dd\te _j \wedge \db\te _j ) \\
 &+ i\sum _{q,j,k}
 |g_{q,j,k}|^{2a}(2 \dd\te _j  +a\frac{\te _j}{g_{q,j,k} } \dd
 g_{q,j,k})\we \overline{(2 \dd\te _j  +a\frac{\te _j}{g_{q,j,k} } \dd
 g_{q,j,k})}.
\endaligned
$$
Since $\te _j $ are smooth the first sum exceeds $-\frac{1}{2c}
\om _Y$ (in the sense of currents) for $c$ small enough. All terms
in the second sum  are  positive. On $V'_j$ we have  $\te _j =1 ,$
the form $\dd\te _j$ vanishes,  and thus we obtain \thetag{1.1}.
\enddemo

We apply Lemma 1 also on  $X$ for
$$
\tilde{\psi } _a = \sum _{j,k} \ze _j ^2  |h_{jk}|^{2a} ,
$$
and fix $c>0$ so that \thetag{1.1} holds for $\psi _a$ and the
corresponding inequality on $X$ is true for $\tilde{\psi } _a$.
Define
$$
\fii _t =c\psi _a\circ \pi +ct\tilde{\psi } _a .\tag{1.2}
$$
By \thetag{1.1}
$$
\fii _t \in PSH(\frac{1}{2}\om _t ) .\tag{1.3}
$$
By the Calabi-Yau theorem [Y] and its non smooth version [K1] it
is no loss of generality to assume that $F$ and $u_t$ are smooth
(provided that a priori estimates do not depend on derivatives of
$u_t$ and $F$) . We treat two cases separately:

\newpage

\bf CASE 1. \rm There exist  $\ep >0$ and $\de >0$ such that for
any $t\in (0,1) $ and any $s\in (0, -s_t)$
$$
\int _{U(t,s, 2\ep )}\om _X^n \geq \de \int _{U(t,s)}\om _X^n ,
$$
where $U(t,s)=\{ u_t < \fii _t  + s_t + s \}, \ s_t =\inf (u_t
-\fii _t )$, $U(t,s, \ep )=U(t,s)\setminus \bar{S^{\ep } } .$

\bigskip
We seek for a uniform bound on $s_t .$ For $s\in (0 ,-s_t]$ define
$$
\Phi _t (s) =\frac{s }{(\int _{U(t,s)}\om _X^n )^{1/n} }.$$ It is
enough to find a bound for $\Phi _t  (s)$ since
$$
\Phi _t (-s_t )=|s_t |/[ \int _{X}\om _X^n ]^{1/n} .
$$
First we show that
$$
\limsup _{s\to 0} \Phi _t (s)\leq C .\tag{1.4}
$$
On the set $X\setminus S^{\ep }$ we have
$$ C(\ep ) \om _X
^{n-m} \we \om ^m \geq \om _X ^n .\tag{1.5}$$
 Hence if $t^{n-m}
F\om _X ^n =f_t \om _t ^n $ then $f_t \leq CF$ on this set.
 Since $u_t $ is smooth and $\fii _t $ is smooth on $X\setminus S^{\ep }$ we have for $\zeta\in X\setminus S^{\ep }$,
  where $u_t -\fii _t$
attains its minimum, $D^2 (u_t -\fii _t  )(\zeta )\geq 0$ and so,
having $\fii _t$ in $PSH (\frac{1}{2} \om _t ),$
$$\aligned
C\om _t^n (\zeta ) &\geq f_t (\zeta )\om _t^n (\zeta )=  (\om _t
+dd^c u_t )^n
(\zeta )  \\
&\geq 2^{1-n} dd^c (u_t -\fii _t ) \wedge \om _t^{n-1} (\zeta )
\geq 2^{1-n} dd^c (u_t -\fii _t ) \wedge \om _t^{n-1} (\zeta
).\endaligned$$
Therefore, interpreting this inequality in local
coordinates diagonalizing $\om _t (\zeta )$ we get
$$
|D^2 (u_t -\fii _t ) (\zeta )|\leq C.
$$

Now, in geodesic coordinates around $\zeta $, consider (for $s$
close to $0 $) the maximal ball $B(\zeta ,r(s))$ contained in
$U(t,s)$ . By the estimate on $D^2 (u_t -\fii _t ) (\zeta )$ and
the Taylor expansion for $z\in \partial B(\zeta , r(s))\cap
\partial U(t,s)$

$$s=u_t (z)-\fii _t (z) -s_t \leq C |z-\zeta |^2 = Cr(s)^2 .$$ So

$$
\int _{U(t,s)}\om _X^n \geq \int _{B(\zeta ,r(s))}\om _X^n\geq
Cr(s)^{2n} \geq Cs^n
$$
and $\Phi _t (s)\leq C$ for $s$ close to zero.

\newpage

Having \thetag{1.4} it is enough to find a uniform bound for $\Phi
_t (a_t )=\max \Phi _t .$ For such $a_t$  we have

$$
\int _{U(t,a_t )}\om _X^n \leq 2^n \int _{U(t, a_t /2  )}\om _X^n
.\tag{1.6}
$$
Using Stokes' theorem we get
$$
\aligned & \int _{U(t, a_t  )}(\om _t +dd^c u_t )\we (\om _t +dd^c
\fii _t )^{n-1} \\ \leq & \int _{U(t, a_t  )}(\om _t +dd^c u_t )^2
\we (\om _t +dd^c \fii _t )^{n-2} \leq ... \\ \leq & \int _{U(t,
a_t )}(\om _t +dd^c u_t )^n .
\endaligned
$$
Integrating by parts, applying \thetag{1.3} and the above
inequality one arrives at

$$\aligned
&\int _{U(t, a_t )}d(u_t  -\fii _t )\we d^c (u_t  -\fii _t )
\we\om
_t^{n-1}  \\
=&\int _{U(t, a_t  )}(s_t +a_t -u_t +\fii _t )dd^c (u_t  -\fii _t )  \we \om _t^{n-1}  \\
\leq & \int _{U(t, a_t  )}(s_t +a_t  -u_t +\fii _t )(\om _t
+dd^c u_t  )\we \om _t^{n-1} \\
\leq  2^n a_t &\int _{U(t, a_t  )}(\om _t +dd^c u_t )\we (\om _t
+dd^c \fii _t )^{n-1} \\
\leq  2^n a_t &\int _{U(t, a_t  )}(\om _t +dd^c u _t )^n \\
=& 2^n a_t t^{n-m} \int _{U(t, a_t )}c_t F\om _X ^n .
\endaligned
$$
From \thetag{1.5} it follows that for any (1,0) form $\gamma$

$$i\gamma \we \bar{\gamma} \we \om _t^{n-1} \geq Ct^{n-m} i\gamma \we
\bar{\gamma}\we \om _X^{n-1}
$$
on $U(t, a_t , \ep  )$.
 By this and the previous inequality   one
obtains
$$\aligned
t^{n-m}&\int _{U(t, a_t , \ep  )}d(u_t  -\fii _t ) \we d^c (u_t
-\fii _t ) \we\om _X^{n-1} \\  \leq C&\int
_{U(t, a_t , \ep )}d(u_t  -\fii _t ) \we d^c (u_t  -\fii _t )  \we\om _t^{n-1} \\
 \leq
C&\int _{U(t, a_t )}d(u_t  -\fii _t ) \we d^c (u_t  -\fii _t )
\we\om _t^{n-1} \leq
 Ca_t t^{n-m}\int _{U(t, a_t)}\om _X ^n .
\endaligned
$$
\newpage

Thus, by the assumption of Case 1,

$$
\int _{U(t, a_t ,\ep  )}d(u_t  -\fii _t ) \we d^c (u_t  -\fii _t )
\we\om _X^{n-1} \leq Ca_t \int _{U(t, a_t)}\om _X ^n \leq  Ca_t
\int _{U(t, a_t  , \ep  )}\om _X ^n .\tag{1.7}
$$
We cover $X\setminus S^{2\ep }$ by a finite number of unit cubes
(in local coordinates) $W_j , j=1,2, ... ,N_0 $ such that $W_j
\cap S^{\ep } =\emptyset $. For one of them, say $W_1$, we have
$$
\int _{U(t, a_t /2 , 2\ep )\cap W_1 } \om _X ^n \geq (1/N_0)\int
_{U(t, a_t /2  ,2\ep )} \om _X ^n .\tag{1.8}
$$
Let $(x_1 , x_2 , ... , x_{2n} )$ denote real coordinates in $W_1
= \{ x: |x_j |<1 \}. $ Set $\Om _1 =W_1 \cap U(t, a_t , \ep ) =W_1
\cap U(t, a_t ),\  \Om _2 =W_1 \cap U(t, a_t /2 , 2\ep  )$ and let
$\pi _j $ denote the projection

$$
\pi _j (x)=(x_1 , ... , x_{j-1} , x_{j+1}, ... , x_{2n}) .
$$
By an isoperimetric inequality from [LW] there is $j$ (and we take
$j=1$) such that
$$
V_{2n-1} (\pi _1 (\Om _2 ))\geq V_{2n} (\Om _2 )^{\frac{2n-1}{2n}
} , \tag{1.9}
$$
where $V_k$ denotes Euclidean volume in $\Bbb R ^k .$ We shall
prove our estimate in two subcases separately.
\bigskip

\bf Case 1A. \rm For $G=\{ y\in \pi _1 (\Om _2 ): \pi ^{-1}_1
(y)\cap \partial \Om _1 \neq\emptyset \} $ we have
$$
V_{2n-1} (G)<\frac{1}{2} V_{2n-1} (\pi _1 (\Om _2 )).
$$
 Observe that for $y\in \pi _1 (\Om _2 )\setminus G$  we have
$\pi ^{-1}_1 (y) \cap W_1 \subset \Om _1$.
 Therefore, by Fubini's theorem
$$
 V_{2n-1} (\pi _1 (\Om _2 )) \leq 2 V_{2n} (\Om _1 ).
$$
Applying the assumption of Case 1, \thetag{1.6 }, \thetag{1.8 },
\thetag{ 1.9} and the last inequality one obtains

$$
V_{2n} (\Om _1 )^{\frac{2n-1}{2n} }\leq C  V_{2n} (\Om _2
)^{\frac{2n-1}{2n} }\leq C V_{2n-1} (\pi _1 (\Om _2 )) \leq C
V_{2n} (\Om _1 ) .
$$
 So $V_{2n} (\Om _1 ) \geq C$ which, via \thetag{1.6}, gives a lower bound for $\int
_{U(t,a_t /2 )}\om _X ^n .$ Therefore
$$
\int _{U(t,a_t /2 )}(-u_t )\om _X ^n \geq C(|s_t |/2 -||\fii _t
||_{\infty } ).
$$
On the other hand, since $\max u_t =0$ there exists (by the
standard Green function argument) a constant $C_0$  such that

\newpage

$$
\int _{X}(-u_t )\om _X ^n \leq C_0
$$
for all $t\in (0,1).$ Combining the last two inequalities we
obtain a uniform bound for $|s_t |$ and further $\Phi _t (a_t
)\leq C.$

\bigskip

\bf Case 1B. \rm Now

$$
V_{2n-1} (G)\geq \frac{1}{2} V_{2n-1} (\pi _1 (\Om _2
)).\tag{1.10}
$$

Let us denote $d(y)=V_1 (G(y)), \ G(y)=\{ x_1\in [-1,1]: (x_1
,y)\in \Om _1 \},$ and observe that for $y\in G$ the set
$G(y)\times \{ y\} $ contains an open  interval in $\Om _1$
joining points from $\Om _2$ and $\partial \Om _1$. Then the
integral of $|\frac{\partial (u_t -\fii _t )}{\partial x_1 }|$
over this interval exceeds $a_t /2 .$ We use this fact and  the
Schwarz inequality to justify the fourth in the following chain of
inequalities. The first one follows from \thetag{1.7 }, the third
one from Fubini's theorem, the fifth from the Schwarz inequality,
the sixth from \thetag{1.9 }, \thetag{1.10}, the seventh from
\thetag{1.6 }, \thetag{1.8 } and the assumption of Case 1 , and
the last one again from the assumption of Case 1.

 $$\aligned
 &a_t \int
_{U(t,a_t ,\ep )}\om _X ^n \geq C\int _{U(t, a_t ,\ep )\cap W_1
}d(u_t -\fii _t )\we d^c (u_t -\fii _t ) \we\om _X^{n-1}\\ &\geq
C\int _{\Om _1 } |\frac{\partial
(u_t -\fii _t )}{\partial x_1}|^2 \, dV_{2n} \\
& \geq C\int _{G }\int _{G(y)}  |\frac{\partial (u_t -\fii _t
)}{\partial x_1 }|^2 \, dx_1 dy \geq Ca_t ^2 \int _{G
}\frac{1}{d(y)}\, dy
\\ & \geq Ca_t ^2 \frac{V_{2n-1}^2 (G)}{V_{2n}(\Om _1 )}
\geq Ca_t ^2 V_{2n}^{1-1/n} (\Om _1 )\geq Ca_t ^2 (\int _{U(t,a_t
, \ep )}\om _X ^n )^{1-1/n} \\ &  \geq Ca_t ^2 (\int _{U(t,a_t
)}\om _X ^n )^{1-1/n}  .
\endaligned
$$
So

$$
\Phi _t (a_t ) =\frac{a_t }{(\int _{U(t, a_t )}\om _X^n )^{1/n}
}\leq C ,
$$
which finishes the proof of Case 1. Note that the bound of $\Phi
_t $ does not depend on derivatives of $u_t$ and so the additional
assumption that $u_t$ be smooth may be dropped.

\bigskip

\bf CASE 2. \rm
 If the assumption of Case 1 is not satisfied then
for any  $\ep >0$ and $\de >0$ there exist  sequences $t(j) \to 0
$ and  $s(j) \in (0, -s_t)$ (those sequences depend on $\ep$ and
$\de$ which is omitted in the notation) such that

\newpage

$$
\int _{U (t(j),s(j), \ep )}\om _X^n \leq \de \int
_{U(t(j),s(j))}\om _X^n . \tag{1.11}$$
From now on we consider
only $t=t(j)$ and $s=s(j)$ which occur in \thetag{1.11}.

By the assumption \thetag{0.2}, \thetag{1.1} and \thetag{1.3} we
have the following estimates of currents on the set $U'_{j_1}\cap
\pi ^{-1} (V'_{j_0 } )\cap S^{\ep }$:

$$\aligned
 &(\om _t +dd^c \fii _t )^n =
[\pi ^{\ast }  ( dd^c c\psi _a +\om _Y )+ t(dd^c c\tilde{\psi } _a
+\om _X )]^n  \\
&\geq t^{n-m} [\pi ^{\ast }  ( dd^c c\psi _a +\om _Y )]^m \we
(dd^c c\tilde{\psi } _a +\om _X )^{n-m} \\
&\geq t^{n-m} \pi ^{\ast } [ \frac{1}{2} \om _Y + ica^2 \sum
_{q,k} |g_{q,j_0 ,k} |^{2(a-1)}
\partial g_{q,j_0 ,k} \we \bar{\partial }\bar{g}_{q,j_0 ,k} ] ^{m}
 \\ &  \we  [ \frac{1}{2} \om _X + ica^2 \sum _{q,k} |h_{q,j_1 ,k}
|^{2(a-1)}
\partial h_{q,j_1 ,k} \we \bar{\partial }\bar{h}_{q,j_1 ,k}  ]^{n-m}\\
&\geq Ct^{n-m}H(\ep )\om _X ^n .
\endaligned\tag{1.12}
$$
Using this and the comparison principle [K2] one obtains

$$\aligned
Ct^{n-m} \int _{U(t,s)} \om _X ^n &\geq t^{n-m} \int _{U(t,s)}
F\om _X ^n = \int _{U(t,s)} (\om _t +dd^c u_t )^n \\ \geq & \int
_{U(t,s)\cap S^{\ep } } (\om _t +dd^c \fii _t )^n  \geq  Ct^{n-m}
H(\ep ) \int _{U(t,s)\cap S^{\ep }  } \om _X ^n .
\endaligned
$$
Combine it with \thetag{1.11} to arrive at

$$\aligned
\int _{U(t,s)\cap S^{\ep }  } \om _X ^n &\geq (1-\de ) \int
_{U(t,s) } \om _X ^n \\ &\geq (1-\de )  CH(\ep ) \int _{U(t,s)\cap
S^{\ep } } \om _X ^n ,
\endaligned
$$
which contradicts the assumption $\lim _{\ep \to 0} H(\ep
 )=\infty $ for small $\ep$. Thus Case 2 never occurs.

\newpage

\proclaim{\bf 2. Proof of Theorem 2} \endproclaim \rm We adopt
notations in the introduction. First observe that

$$\omega(t,\cdot)=e^{-t} \omega_0 + (1-e^{-t})\pi^*\omega_Y + i\partial\overline{\partial}
\varphi(t,\cdot).$$
Denote by $\omega_t$ the K\"ahler form $e^{-t}
\omega_0 + (1-e^{-t})\pi^*\omega_Y$. Adding appropriate constants
to $\varphi(t,\cdot)$, one can show

$$\frac{\partial\varphi}{\partial t} = \log\frac{e^{(n-m)t}
(\omega_t+ i\partial\bar\partial \varphi)^n}{\Omega
}- \varphi, ~~~~\varphi(0,\cdot)=0, $$
 where $\Omega$ is a volume form on $X$ determined by
$\text{ Ric }(\Omega)=-\pi^*\omega_Y$. Then there is a constant
$c$ such that

$$e^{(n-m)t}\omega_t^n = \sum_{j=1}^m {n \choose j} (1- e^{-t} )^j e^{(j-m)t}
\pi^*\omega_Y^j \wedge \omega_0^{n-j} \le c \Omega.$$ Therefore,
by the Maximum principle, we have $\varphi \le {\max}\{\log
c,0\}$.
Differentiating the above equation on $t$, we get

$$\frac{\partial}{\partial t}\left( \frac{\partial \varphi}{\partial t} \right) =
\Delta _t \frac{\partial \varphi}{\partial t} - \frac{\partial
\varphi}{\partial t} + \text{ tr } _t ( e^{-t}(\pi^*\omega_Y -
\omega_0)) + (n-m).
$$
Again by using the Maximum principle, one can show that
$\frac{\partial \varphi}{\partial t}$ is uniformly bounded from
above.

Set $s=e^{-t}/(1-e^{-t})$, we then have
$$(s \omega_0 + \pi^*\omega_Y + i\partial\bar\partial \psi)^n = s^{n-m} F \omega_0^n,$$
where $\psi = (1-e^{-t})^{-1} \varphi$ and

$$F = (1-e^{-t})^{-m} e^{\varphi + \frac{\partial \varphi}{\partial t}} \frac{\Omega}{\omega_0^n}.$$
Clearly, $F$ is uniformly bounded, so we can apply Theorem 1 to
conclude that $\psi$ is uniformly bounded, so does $\varphi$.

The rest of the proof is exactly the same as that in [ST1] since
all other estimates there are dimension free (also see [ST2]).
This finishes the proof of Theorem 2.

\newpage

\proclaim{3. Examples}\endproclaim\rm In this section, we verify
the assumptions \thetag{0.2} for two examples. We will adopt
notations in the introduction. First we assume that $\pi: X\mapsto
Y$ is a holomorphic fibration with $\text{ dim } Y = 1$, that is,
all fibers are hypersurfaces in $X$. In this case, any point of
$Y$ has a neighborhood which can be identified with an open subset
$U$ in the complex plane and $\pi$ is given by a holomorphic
function $h$. By \thetag{0.3}, we only need to show for some $a\in
(0,1)$

$$|h-h(x)|^{2(a-1)} i\partial h \wedge \bar{\partial } h \wedge
\omega_X^{n-1} \ge H(\epsilon) \omega^n_X.$$ in a small
neighborhood of $x$, which may depend on $\epsilon$, where $x$ is
any fixed point in $\pi^{-1}(U)$ and $H(\epsilon)$ was given in
\thetag{0.2}. By a result of \L ojasiewicz [Lo], there are
constants $\theta\in (0,1/2)$ and $\sigma_1$ such that for any
$\xi$ with $d(x,\xi)\le \sigma_1$,

$$|\nabla h|(\xi)\ge |h(\xi)-h(x)|^{1-\theta} .\tag{3.1}$$
Now choose $a < \theta$, then \thetag{3.1} holds for $h$. Therefore
 the assumptions in \thetag{0.2} hold for any
1-dimensional base.
\bigskip

Before we discuss the second example we state a simple lemma
regarding to the assumptions \thetag{0.2}.

\proclaim{Lemma 2} Let $\pi:X\mapsto Y$ be the product of a
holomorphic fibration $\pi_1: X_1\mapsto Y_1$ and a complex
manifold $U$, that is, $X=X_1\times U$, $Y=Y_1\times U$ and
$\pi=\pi_1\times {\text{Id}_U}$. If the fibration
$\pi_1:X_1\mapsto Y_1$ satisfies \thetag{0.2}, then so does $\pi:
X\mapsto Y$.
\endproclaim
It follows directly from \thetag{0.2} since $\pi$ is
non-degenerate along $U$-directions.

\bigskip

Now we consider the second example. Let $\pi:X\mapsto Y$ be a
generic fibration over a complex surface $Y$ with elliptic curves
as fibers. Let $A$ be the set of $y\in Y$ such that $\pi^{-1}(y)$
is a singular elliptic curve. Since the fibration is generic, $A$
is a divisor in $Y$ and each singular fiber has either a node or
an ordinary cusp.

\proclaim{Proposition 3} Such an elliptic fibration $\pi: X\mapsto
Y$ over a surface satisfies \thetag{0.2}.
\endproclaim

The rest of this section is devoted to the proof of this
proposition.

Set $A_0$ to be the subset of $A$ over which singular fibers have
nodes and $A_1$ to be the subset of those points $y\in A$ such
that $\pi^{-1}(y)$ has a cusp. Clearly, $A_1$ consists of finitely
many isolated points. As usual, let $S$ be the set of singular
points in the singular fibers of $\pi: X\mapsto Y$. Then
$S_0=\pi^{-1}(A_0)$ and $S_1=\pi^{-1}(A_1)$. By Lemma 2, clearly,
\thetag{0.2} holds for any point in $S_0$. So we only need to
check \thetag{0.2} for any $p\in S_1$.

Fix each such a $p$, we can find local coordinates $x,y,z$ of $X$
near $p$ and local coordinates $s,t$ of $Y$ near $\pi(p)$
satisfying: $\pi (x,y, z )= (z ,  y^2 -x^3 +zx )=(s,t)\in Y$.
Furthermore, we may assume

$$ \om _X =i( dx\we d\bar{x} +dy\we d\bar{y} + dz\we
d\bar{z})$$

\newpage
and

$$ \om _Y =i( ds\we d\bar{s} +dt\we d\bar{t}).$$
Then we have
$$
\om = \pi ^{\ast} \om _Y =i (dz\we d\bar{z} +[(z-3x^2 )dx +2y
dy]\we \overline{[(z-3x^2 )dx +2y dy]} )
$$ and
$$
\om ^2 =-2 dz\we d\bar{z} \we [(z-3x^2 )dx +2y dy]\we
\overline{[(z-3x^2 )dx +2y dy]}.
$$
So $$ \om ^2 \we \om _X = (-2i) (|z-3x^2 |^2 +4|y|^2 ) dx\we
d\bar{x}\we dy\we d\bar{y} \we dz\we d\bar{z} .
$$
The singular set is

$$
S=  \{ y=0 \}  \cap \{ z=3x^2 \}.
$$
Since $\pi (x,0,3x^2 )= (3x^2 ,2x^3 )$  the image of $S$ lies in
$$
 A= \{ 4s^3 =27t^2 \} .$$
As the generator for $A$ close to origin we take $g(s,t)=4s^3
-27t^2$, and as generators for $S$ we take
$$
h_1 (x,y,z)= y, \ \ \ h_2 (x,y,z) =z-3x^2 .
$$
Define, for positive $a<1/100$
$$\aligned
\ga _1 &= \pi ^{\ast} (i |g|^{a-2} \partial g \we \bar{\partial
}\bar{g} \we  \om _Y ) \\
=& \pi ^{\ast} (36i |4s^3 -27 t^2|^{a-2}(4|s|^4 +81 |t|^2 ) ds\we
d\bar{s} \we dt\we d\bar{t} ). \endaligned
$$ and

$$\aligned
\ga _2 &=i( |h_1 |^{a-2} \partial h_1 \we \bar{\partial }\bar{h}_1
+  |h_2 |^{a-2} \partial h_2 \we \bar{\partial }\bar{h}_2 ) +\om
_X
\\
&\geq i|y|^{a-2}dy\we d\bar{y} +i |z-3x^2 |^{a-2} (dz-6xdx)\we
\overline{(dz-6xdx)} +idx\we d\bar{x}.\endaligned
$$
Wedging $\ga _1$ and $\ga _2$ we can forget about differentials in
$\ga _2$ containing $dz$ or $ d\bar{z}$. Thus we seek for the
lower bound for $ \ga_1 \we \ga _3$ , where
$$
\ga _3 = i|y|^{a-2}dy\we d\bar{y} +i (36 |x|^2|z-3x^2 |^{a-2}
+1)dx\we d\bar{x} .$$ We have

$$\aligned
\ga_1 \we \ga _3 &\geq \pi ^{\ast} [(i |4s^3 -27 t^2|^{a-2}(4|s|^4
+81 |t|^2 )] dz\we d\bar{z} \we \pi ^{\ast} [idt\we d\bar{t}]\we
\ga _3 \\ &=  \pi ^{\ast} [(i |4s^3 -27 t^2|^{a-2}(4|s|^4 +81
|t|^2 )] dz\we d\bar{z} \\ &\we i[(z-3x^2 )dx +2y dy]\we
\overline{[(z-3x^2 )dx +2y dy]}\we \ga _3
\\ &\geq \pi ^{\ast} [(i |4s^3 -27 t^2|^{a-2}(4|s|^4 +81 |t|^2 )]
\\ &\times
[|z-3x^2 |^2 |y|^{a-2} +4|y|^2 (36 |x|^2 |z-3x^2 |^{a-2}
+1)](1/6)\om _X ^3
\endaligned$$

\newpage
Let us denote
$$\aligned
f_1 (x,y,z)&= \pi ^{\ast} ( |4s^3 -27 t^2|), \\
f_2 (x,y,z)&= \pi ^{\ast} (4|s|^4 +81 |t|^2 ),
\endaligned$$
and writing $w=z-3x^2$:
$$\aligned
f_3 (x,y,z)&= |z-3x^2 |^2 |y|^{a-2} +4|y|^2 (36 |x|^2 |z-3x^2
|^{a-2} +1) \\
&=   |w|^2 |y|^{a-2} +4|y|^2 (36 |x|^2 |w |^{a-2} +1).
\endaligned $$
To check the hypothesis of the theorem we need to show that
$$
\lim f_1 ^{a-2} f_2 f_3 =\infty
$$
when $(x,y,z)$ tends to the origin. Since
$$
f_1 \leq 100 \max (|s|^3 , |t|^2 )$$ we conclude that
$$
f_1 ^{-4/3 -a} f_2 \to\infty .
$$
To finish the verification we need to prove that $ f_1 ^{-2/3 +2a}
f_3$ is bounded away from zero. It is enough to check it for
points satisfying $f_3 <1$. So we assume that
$$
|w|^2 <|y|^{2-a} \ \ \ \text{and so} \ \ |w|^{a-2} \geq |y|^{-5/3}
. \tag{3.2}
$$
We have
$$
f_1 =|w^2 (4w +9x^2 )- 27y^2 (4x^3 +2wx +y^2 ) |.
$$
Hence
$$
f_1 \leq M:= 500 \max (|w|^3 , |w^2 x^2 |, |y^2 x^3 |, |y|^4 ).
$$
Accordingly consider four cases:

1) $M=500 |w|^3 $ (then $|w^3 | \geq |y^4 | $). So
$$
500 f_1 ^{-2/3 +2a} f_3 \geq |w|^{6a-2} |w|^2 |y|^{a-2} \geq
|w|^{6a-1} |y|^{a-2/3} \to\infty .
$$

2) $M=500 |xw|^2 $  (then $|w | \leq |x^2 | $). So (using $|y^2
\geq |w|^{2+2a}$ which follows from \thetag{3.2})
$$
500 f_1 ^{-2/3 +2a} f_3 \geq |wx|^{4a-4/3} |xy|^2  |w|^{a-2} \geq
|x|^{4a-1/3}  |w|^{7a-5/6}\to\infty .
$$

3) $M=500 |x^3 y^2 | $  (then $|y^2 | \leq |x^3 | $). So (see
point 1))
$$\aligned
500 f_1 ^{-2/3 +2a} f_3 &\geq |y|^{4a-4/3} |x|^{6a-2} |xy|^2
|w|^{a-2} \geq |y|^{4a-4/3}  |y|^{2+4a} |w|^{a-2}\\ & \geq
|y|^{-1+8a}\to \infty . \endaligned$$

4) $M=500 |y^4 | $. So (see \thetag{1})$$ 500 f_1 ^{-2/3 +2a} f_3
\geq |y|^{8a-8/3} |y|^2 \to\infty .$$

Thus Proposition 3 is proved.

\newpage

\bigskip

\bf\centerline{References}\bigskip \rm

\bigskip

\noindent\item {[DP] }{J.-P. Demailly, M. Paun, \it Numerical
characterization of the K\"ahler cone of a compact K\"ahler
manifold, \rm Ann. Math., \bf 159 \rm (2004), 1247-1274.}

\noindent\item{[K1] }{S. Ko\l odziej, \it The complex \MA
equation, \rm Acta Math. \bf 180 \rm (1998), 69-117.}

\noindent\item {[K2] }{S. Ko\l odziej, \it   Stability of
solutions to the complex Monge-Amp\`ere on compact K\"ahler
manifolds,  \rm Indiana U. Math. J. \bf 52 \rm (2003), 667-686.}

\noindent\item {[Lo] }{S. \L ojasiewicz, \it Emsembles
semi-analtiques, \rm I.H.E.S. notes (1965).}

\noindent\item {[LW] }{L.H.  Loomis, H. Whitney, \it An inequality
related to the isoperimetric inequality, \rm Bull. Amer. Math.
Soc. \bf 55 \rm (1949), 961-962.}

\noindent\item {[ST1] }{J. Song and G. Tian, \it The
K\"ahler-Ricci flow on surfaces of positive Kodaira dimension, \rm
To appear in Inventiones Math.}

\noindent\item {[ST2] }{J. Song and G. Tian, \it Generalized
K\"ahler-Einstein metrics, \rm in preparation.}

\noindent\item {[TZ] }{G. Tian and Z. Zhang, \it On the
K\"ahler-Ricci flow on projective manifolds of general type, \rm
Chinese Ann. Math. B \bf 27 (2) \rm (2006), \rm 179-192.}

\noindent\item {[Y] }{S.-T. Yau, \it On the Ricci curvature of a
compact K\"ahler manifold and  the complex Monge-Amp\`ere
equation, \rm Comm. Pure and Appl. Math. \bf 31 \rm (1978),
339-411.}

\bigskip

\bigskip

\bigskip

\bye